\documentclass[12pt]{article}
\usepackage{latexsym}
\begin{document}
\title{Topology of Solutions of the Liouville Equation\\}
\author{ W{\l}odzimierz Piechocki \vspace{3mm}\thanks{Presented 
at the XVI Workshop
on Geometric Methods in Physics, Bia{\l}owie\.{z}a, Poland, 
June 30 - July 6, 1997}\\ \small Field Theory Group,
So{\l}tan Institute for Nuclear Studies,\\\small Ho\.{z}a 69, 
00-681 Warsaw, Poland. E-mail: piech@fuw.edu.pl \\}
\maketitle
\begin{abstract}
Suggestions concerning the generalization of the geometric 
quantization to the case of nonlinear field theories are given. 
Results for the Liouville field theory are presented.
\end{abstract}
\pagebreak

Significant part of cosmic rays is connected with the existence of
black holes predicted by classical theory of gravitation. However, 
correct description of these data can be done only in the framework 
of quantum theory. Elementary particles produced in laboratories  
are well described by the standard model. This model, however, is 
a phenomenological one; it has nearly 20 free parameters. 

To describe satisfactory both celestial and terrestrial data we should 
quantize, in a rigorous way, the Einstein general relativity and the 
Yang-Mills theory. Both theories are \emph{nonlinear field} theories. 
The problem is that such theories have \emph{infinitely} many degrees 
of freedom and the space of solutions to the field equations 
\emph{is not} a vector space.

We have a few methods of quantization. Geometric quantization~\cite{1} 
is a method that has sound mathematical foundation and properly 
describes simple systems with \emph{finitely} many degrees of freedom. 
It is not clear, however, if we can generalize this method to the case 
of nonlinear field theories.

Both general relativity and Yang-Mills theories are rather \emph{
complex} theories. It is perhaps reasonable to apply geometric 
quantization first to \emph{simple} integrable systems to establish 
the framework. An example of such a system is the 2-dim
Liouville theory~\cite{2}. The initial value problem for the Liouville 
equation reads:
\begin{equation}
\left( \partial^2_{t} - \partial^2_{x} \right) \varphi (t,x) + 
\frac{m^2}{2} \exp \varphi (t,x) = 0,\;\;\;\; m > 0
\end{equation}
\begin{equation}
\left\{ \begin{array}{lll}\varphi (0,x)&=&\phi (x)  \\
\varphi_t (0,x) &=&\pi (x)
\end{array} 
\right.
\end{equation}
Equation (1) appears, e.g., in the context of the $SU(2)$ Yang-Mills 
theory~\cite{3}. Fields which minimize locally the Yang-Mills finite 
action in 4-dim Euclidean space have to satisfy the equation 
$F_{\mu \nu} = \widetilde{F}_{\mu \nu}\;$ (where $\widetilde
{F}_{\mu \nu}
:= \frac{1}{2}\epsilon_{\mu \nu \lambda \gamma} F^{\lambda \gamma}$ ).
Solutions to this equation invariant under 3-dim rotations combined 
with gauge transformations satisfy Eq.(1). Complete solution of the 
Cauchy problem for the Liouville equation was obtained in late 70-ties
~\cite{4}. One can easily check~\cite{4} that Eq.(1) is satisfied by 
the function:
\begin{equation}
\varphi (t,x) := - \log \left[ \frac{m^2}{16}F^2(t,x)\right],
\end{equation}
where
\begin{equation}
F(t,x) := \chi_1(x+t)\psi(x-t)+\chi_2(x+t)\psi_1(x-t),
\end{equation}
and where $ \chi_k, \psi_k \;\;(k=1,2)$ are any functions that fulfill 
two requirements:
\begin{equation}
\left \{ \begin{array}{lll}\chi_1 \chi_2^{\prime} - \chi_2 \chi_1^
{\prime}&=&1\\ \psi_1 \psi_2^{\prime}-\psi_2 \psi_1^{\prime}&=&1
\end{array}
\right.
\end{equation}
The solution (3) can be singular if there exists $(t_0,x_0)\in R^2$ 
such that $F(t_0,x_0)=0$. Applying the implicite mapping theorem to 
Eqs. (3-5) one can prove~\cite{5,4} the lemma:
\[
 \left( \exists \;(t_0,x_0) \in R^2\;:\;F(t_0,x_0) = 0 \right)
\longrightarrow \left( \begin{array}{l}\exists\;\;
R\ni t \rightarrow x(t) \in R \;\;\mbox{such that}\\
x(t_0)=x_0 \;\;\mbox{and}\;\;F(t,x(t))=0 
\end{array} \right)  \]
The domain of the mapping $x(\cdot)$ is the entire set of real numbers.

Solutions that satisfy the initial value data (2) specify $\psi_k$ 
and $\chi_k$\\( see~\cite{4}):
\begin{equation}
\left\{ \begin{array}{l}\psi_k^{\prime \prime}= u \psi_k \\
\chi_k^{\prime \prime}=w\chi_k 
\end{array}
\right.    \;\;\;\;\mbox{for}\;\;k=1,2
\end{equation}
where
\begin{equation}
\left\{ \begin{array}{l}
u:= \frac{1}{16}\left[ (\phi^{\prime}-\pi)^2-4(\phi^{\prime}-\pi)
^{\prime}
+ m^2 \exp \phi \right] \\
w:= \frac{1}{16}\left[ (\phi^{\prime}+\pi)^2 - 4(\phi^{\prime}+\pi)
^{\prime}+m^2\exp\phi \right] \end{array}  \right.
\end{equation}

Corollary resulting from the lemma is that
\begin{equation}
\left( (\phi,\pi) \in C^{\infty}(R)\times C^{\infty}(R) \right)
\Longrightarrow \left( \varphi \in C^{\infty}(R^2,R) \right)
\end{equation}

From now on, we only consider the smooth solutions.

The space of initial data $\mathcal{F}:= C^{\infty}(R)\times C^{\infty}
(R)$ is a Fr\'{e}chet space. However, the set of solutions
\[ \mathcal{M}:= \left\{\varphi \in C^{\infty}(R^2,R)\;\mid \; \varphi 
\;\mbox{satisfy Eqs.(1) and (2)} \right\}\;\subset C^{\infty}(R^2,R) \]
cannot be a Fr\'{e}chet space, since it is not a vector space.

It turns out that $\mathcal{F}$ and $\mathcal{M}$ are \emph
{homeomorphic}. The proof is elementary but lenghty. Details are given 
in Ref. [5]. Consecutive steps of the proof are roughly the following:
\begin{enumerate}
\item One divides the mapping $S:\mathcal{F} \rightarrow \mathcal{M}$ 
into a few mappings:
\[ \mathcal{F} \ni (\phi,\pi) \stackrel{S_1}{\rightarrow} (u,w)
\stackrel{S_2}{\rightarrow} (\psi,\chi)\stackrel{S_3}{\rightarrow}
\ldots \stackrel{S_n} {\rightarrow}\varphi \in \mathcal{M} \]
and proves that $S_1,S_2,\ldots ,S_n$ are continuous. Then, 
$S:=S_n \cdot \ldots  \cdot S_2
\cdot S_1$ is continuous as a composition of continuous mappings.
\item The inverse mapping is continuous as it is defined by:
\[ S^{-1}:\mathcal{M} \ni \varphi \Longrightarrow S^{-1}(\varphi)
:=\left( \varphi(0,\cdot),\varphi_t(0,\cdot) \right) \in \mathcal{F}\]
\item Since we consider the solutions of the Cauchy problem we get 
\[ S\cdot S^{-1} = \mathbf{1}_{\mathcal{M}}\;\;\;\; \mbox{and}\;\;\;\;
S^{-1}\cdot S
=\mathbf{1}_{\mathcal{F}} \]
\end{enumerate}
This completes the proof.

The geometric quantization procedure assumes that the phase space 
of a given classical system has the structure of a \emph{manifold}. 
In the case of a simple system with n-degrees of freedom the phase 
space is a manifold modelled on $R^{2n}$. 
In the case of a field theory we would 
like to have an object that we could call a manifold modelled on a 
Fr\'{e}chet space. The homeomorphism of $\mathcal{F}$ and $\mathcal{M}$
 gives $\mathcal{M}$ the structure of a topological manifold modelled 
on $\mathcal{F}$. There is only a single chart $(\mathcal{M},S^{-1})$. 

The starting point in the geometric quantization procedure of a 
mechanical system is to express the evolution of the system on the 
phase space in terms of \emph{symplectic}  geometry. Can one follow 
this procedure in the case of the Liouville field theory?  We shall 
try to answer this question in the near future~\cite{6}.

\section*{Acknowledgements}
I would like to thank the organizers of the Workshop for a pleasant 
and stimulating atmosphere. I also wish to thank Professor M. Flato 
for inspiration.

\end{document}